\documentclass[11pt]{amsart}
\usepackage{tkz-euclide}
\usepackage{graphicx, color}
\usepackage{amscd}
\usepackage{amsmath}
\usepackage{amsfonts}
\usepackage{amssymb}
\usepackage{mathrsfs}
\newtheorem{theorem}{Theorem}[section]

\newtheorem{ex}{Example}[section]

\newtheorem{remark}{Remark}[section]

\date{\today}

\numberwithin{equation}{section}
\newcommand{\R}{\mathbb R}

\newcommand{\NN}{\mathbb N}

\def\Fix{\mathop{\rm Fix}}
\begin{document}
\title[On nonlinear Schr\"{o}dinger equations on the hyperbolic space]{On nonlinear Schr\"{o}dinger equations \\
on the hyperbolic space}
\author[M. Cencelj]{M. Cencelj}
\address[M. Cencelj]{Faculty of Education and Faculty Mathematics and Physics, University of Ljubljana \& Institute of Mathematics, Physics and Mechanics, Ljubljana, Slovenia}
\email{matija.cencelj@fmf.uni-lj.si}
\author[I. Farag\'{o}]{I. Farag\'{o}}
\address[I. Farag\'o]{Department of Differential Equations, Budapest University of Technology and Economics, Budapest, Hungary \& Department of Applied Analysis and Computational Mathematics, E$\ddot{\mbox{o}}$tv$\ddot{\mbox{o}}$s Lor\'{a}nd University, Budapest, Hungary \& MTA-ELTE NumNet Research Group, Budapest, Hungary}
\email{\tt faragois@cs.elte.hu}
\author[R. Horv\'{a}th]{R. Horv\'{a}th}
\address[R. Horv\'{a}th]{
Department of Analysis, Budapest University of Technology and Economics, Budapest, Hungary
\& MTA-ELTE NumNet Research Group, Budapest, Hungary}
\email{\tt rhorvath@math.bme.hu}
\author[D. Repov\v{s}]{D. Repov\v{s}}
\address[D. Repov\v{s}]{Faculty of Education and Faculty Mathematics and Physics, University of Ljubljana \& Institute of Mathematics, Physics and Mechanics, Ljubljana, Slovenia}
\email{\tt dusan.repovs@guest.arnes.si}
\keywords{Elliptic problem on Hadamard manifolds,
Laplace-Beltrami operator,
Poincar\'{e} ball model,
variational methods, energy functional,
growth condition,
Palais principle of symmetric criticality,
Kirchhoff-type problems,
Sobolev space,
 multiple solutions.\\
\phantom{aa} {\sl 2010 Mathematics Subject Classification}:
Primary: 35J60, 49K10, 58J32; Secondary: 35A01, 35R01.}
\begin{abstract}
We study  existence of weak solutions for certain classes of
nonlinear Schr\"{o}dinger equations on the Poincar\'{e} ball model $\mathbb{B}^N$, $N\geq 3$.
By using the Palais principle of symmetric criticality and suitable group theoretical arguments, we establish the existence of a nontrivial (weak) solution.
\end{abstract}
\maketitle

\section{Introduction}\label{sec:Introduction}

This paper was motivated by a large number of applications of the eigenvalues problem for the Laplace-Beltrami operator in the hyperbolic framework and in particular, by recent important work \cite{Molica,MoRa,MoSec}. We study the following elliptic problem
 \begin{equation}\label{BWP6}
 \begin{aligned}
&-\Delta_{H}u=\lambda \alpha(\sigma)f(u) \ \ {\rm on}\,\,{\mathbb{B}^N},
& \ \ u\in H^{1,2}(\mathbb{B}^N),
\end{aligned}
\end{equation}
\noindent on  the Poincar\'{e} ball model $\mathbb{B}^N$.
Here,
$\Delta_{H}$ is the Laplace-Beltrami operator,
$\lambda >0$ is a real parameter,
 $\alpha\in L^1(\mathbb{B}^N)\cap L^\infty(\mathbb{B}^N)$ is a nonnegative nontrivial radially symmetric potential,
$N\geq 3$,
 and
$f:\mathbb{R}\rightarrow \mathbb{R}$ is a continuous function satisfying the following growth condition
\begin{equation}\label{Growth}
\alpha_f=\sup_{t\in \R}\frac{|f(t)|}{1+|t|^{q-1}}<\infty,
\end{equation}
where $q\in \left[2,2^*\right]$ and $2^*=2N/(N-2)$ denotes the critical Sobolev exponent.

Problem \eqref{BWP6} is an important generalization of the most widely investigated elliptic problems with subcritical
nonlinearities which arise naturally in various areas of mathematics.
For instance, an important incentive to study Kirchhoff-type problems comes from recent publications \cite{p22, p23, fisc1, fisc2, mingqi0, mingqi1, mingqi2, MolicaPucciBook, p21} in which Kirchhoff equations on $\mathbb{B}^N$ have been proposed as an interesting open problem (see also \cite{H2,H3,H1,H4,MolicaNA} for related topics).

Since $\mathbb{B}^N$ is an important model of a Hadamard manifold (i.e.\ a complete, simply connected Riemannian manifold with nonpositive sectional curvature), our approach can be used (as we plan to do in our forthcoming paper)  to study  existence of multiple solutions of elliptic problems on Hadamard manifolds in the presence of a compact topological group action.

Given $\sigma\in \mathbb{B}^N$, let $T_\sigma (\mathbb{B}^N)$ denote the tangent space and $\left\langle \cdot, \cdot \right\rangle_\sigma$  the related inner product. We investigate weak solutions of
problem \eqref{BWP6} , i.e.
for  functions $u\in H^{1,2}(\mathbb{B}^N)$ such that for every $\varphi\in H^{1,2}(\mathbb{B}^N)$, the following is satisfied
$$
\int_{\mathbb{B}^N}\left\langle \nabla_H u(\sigma),\nabla_H \varphi(\sigma) \right\rangle_\sigma d\mu =\lambda\int_{\mathbb{B}^N}\alpha(\sigma)f(u(\sigma))\varphi(\sigma)d\mu,
$$
where
$d\mu$ is for the Riemmanian volume element on $\mathbb{B}^N$, and
$
\nabla_H=\left(\frac{(1-|\sigma|^2)}{2}\right)^2\nabla
$
is the covariant gradient (here, $|\cdot|$ and $\nabla$ denote the Euclidean distance and the gradient in $\mathbb{R}^N$, respectively).
 Let $SO(N)$ be the special orthogonal group, $N\geq 3$.
 
We are now ready to state the main result of this paper.
\begin{theorem}\label{MainTheorem}
 Suppose that
 $f:{\R}\to{\R}$
 is a a continuous function such that
 \begin{enumerate}
 \item[(a)] $f$ satisfies the growth condition \eqref{Growth}   for some $q\in (2,2^*)$,
 \item[(b)] $f$ satisfies the asymptotic condition
\begin{equation}\label{azero}
\lim_{t\rightarrow0^+} \frac{\displaystyle\int_0^{t}f(s)ds}{t^2}=+\infty,
\end{equation}
\item[(c)]  $\alpha\in L^1(\mathbb{B}^N)\cap L^\infty(\mathbb{B}^N)\setminus\{0\}$ is a nonnegative radially symmetric map with respect to the origin $\sigma_0\in \mathbb{B}^N$.
\end{enumerate}
Then there exists $\lambda^{\star}>0$
such that, for every $\lambda\in (0,\lambda^{\star})$,
problem \eqref{BWP6} admits a $SO(N)$-invariant weak solution $u_\lambda\in H^{1,2}(\mathbb{B}^N)$ whose norm converges to zero as $\lambda$ goes to zero.
\end{theorem}

We shall find solutions of problem \eqref{BWP6} as critical points of the following energy functional
\begin{equation}\label{Jlambda}
J_\lambda(u)=\dfrac{1}{2}\int_{\mathbb{B}^N}  |\nabla_H u(\sigma)|^2 d\mu-\lambda\int_{\mathbb{B}^N} \alpha(\sigma)\left(\int_0^{u(\sigma)}f(t)dt\right)d\mu
\end{equation}
defined on the Sobolev space $H^{1,2}(\mathbb{B}^N)$. In fact, we shall prove Theorem \ref{MainTheorem} by using variational methods (see  \cite{brezis} as a general reference for this topics) by means of a local minimum result for differentiable functionals, and the Palais principle of symmetric criticality (see Theorems \ref{BMB} and \ref{PalaisTh} below, respectively).

\begin{remark}
Note that condition \eqref{azero} in
Theorem \ref{MainTheorem}
has been used before - in order  to study existence and multiplicity results for certain classes of elliptic problems on bounded domains (see e.g.,  
\cite{Bon}-\cite{BoPi}, 
 \cite{kagi2}-\cite{kagi1},
  and 
  \cite{GMBSe}).
\end{remark}
The noncompact hyperbolic setting presents additional difficulties with respect to the cited work and appropriate geometric and algebraic tools are needed for the proof. A key tool is a detailed analysis of the energy level of $J_\lambda$ on the Sobolev space
$$
H_{SO(N)}^{1,2}(\mathbb{B}^N)=\left\lbrace u\in H^{1,2}(\mathbb{B}^N)\mid g\circledast_{SO(N)}u=u, \mbox{ for every }g\in SO(N)\right\rbrace
$$
of $SO(N)$-invariant functions (see Section \ref{Section3}). 

A simple prototype of a function in $H_{SO(N)}^{1,2}(\mathbb{B}^N)$, depending on parameters $0<r<\rho$, defined by setting for every
 $\sigma\in \mathbb{B}^N$,
\begin{equation}  
w^{1/2}_{\rho,r}(\sigma)=\left\{
\begin{array}{ll}
1 & \mbox{ if $\sigma \in A_{1/2 r}^{\rho}$}\\[5pt]
0 & \mbox{ if $\sigma \in \mathbb{B}^N \setminus A_r^{\rho}$} \\[5pt]
 \displaystyle\frac{2}{r}\left({r- \left|\displaystyle\log\left(\frac{1+|\sigma|}{1-|\sigma|}\right)-\rho\right|}\right)
& \mbox{ if $\sigma \in A_{r}^{\rho}\setminus A_{1/2 r}^{\rho}$}.
\end{array}
\right.
\end{equation}
has the support contained in the annuls $A_r^{\rho}$ of $\mathbb{B}^N$ (see \cite{Molica} for more details).

We conclude the introduction by describing the  structure of the paper. In Section \ref{Section2} we shall collect the necessary
notations, definitions and facts.
In Section \ref{Section3} we shall
present
 a compactness argument, based on the action of a suitable subgroup of the group of isometries of $\mathbb{B}^N$. In
 Section \ref{Section4} we shall prove the main result (Theorem \ref{MainTheorem}). Finally, in Section \ref{Section5} we shall give an example.

Some of the abstract tools used in this paper can be found in \cite{PRR}.
For eigenvalue problems on the hyperbolic space we refer the reader to \cite{BoKl, partenza, 9, MolicaNA, punzo, 15,16}.

\section{The Abstract framework}\label{Section2}

Let $
\mathbb{B}^N=\{\sigma=(x_1, x_2,\dots, x_N)\in \mathbb{R}^N\mid \ |\sigma|<1\},
$
be
equipped  by the Riemannian  metric
$
g_{ij}={4}{(1-|\sigma|^2)^{-2}}\delta_{ij},
$
where $\sigma\in \mathbb{B}^N,  i,j\in\{1,...,N\},$ and  $|\cdot|$ and $\delta_{ij}$ denote the Euclidean distance and the Kronecker delta symbol,  respectively.
For every $i,j\in\{1,\dots,N\}$, let
$
g^{ij}=(g_{ij})^{-1}$ and
 $\textrm{g}={\rm det}(g_{ij}).
$
We locally define the Laplace-Beltrami operator $\Delta_H$   by
$$
\Delta_H=\rm{g}^{-{1/2}}\sum_{i=1}^{N}\frac{\partial}{\partial x_i}\left(\rm{g}^{1/2}\sum_{j=1}^{N}g^{ij}\frac{\partial}{\partial x_j}\right).
$$
The following is a more convenient form
$$
\Delta_H=\frac{(1-|\sigma|^2)^2}{4}\sum_{i=1}^{N}\frac{\partial^2}{\partial x_i^2}+
\frac{(N-2)(1-|\sigma|^2)}{2}\sum_{i=1}^{N}x_i\frac{\partial}{\partial x_i},
$$
when we consider the Riemannian volume element in $\mathbb{B}^N$
\begin{equation}
\label{dmu}
d\mu=\sqrt{\textrm{g}}dx=2^N(1-|\sigma|^2)^{-N}dx,
\end{equation}
where $dx$ denotes the Lebesgue measure on $\R^N$.
Finally, let
\begin{equation}
\label{distance0}
d_H(\sigma)=d_H(\sigma,\sigma_0)=2\int_0^{|\sigma|}\frac{dt}{1-t^2}=\log\left(\frac{1+|\sigma|}{1-|\sigma|}\right)
\end{equation}
be the geodesic distance of  $\sigma\in \mathbb{B}^N$ from the origin $\sigma_0\in \mathbb{B}^N$.
Let $(\varrho,\theta)$ denote the polar geodesic coordinates of a point in $\mathbb{B}^N\setminus\{0\}$. We see that
$
ds^2=d\varrho^2+(\sinh \varrho)^2d\theta,
$
 on $\mathbb{B}^N\setminus\{0\}$,
and
$$
\Delta_H=\frac{\partial^2}{\partial \varrho^2}+(N-1)\coth \varrho\frac{\partial}{\partial \varrho}+\frac{\Delta_\theta}{(\sinh\varrho)^2},
$$
where $\Delta_\theta$ the Laplace-Beltrami operator on the sphere $\mathbb{S}^{N-1}\subset \R^N$.
Invoking \eqref{distance0}, we define the distance on $\mathbb{B}^N$ by
\begin{equation*}
d_{H}(\sigma_1,\sigma_2)
=\operatorname{Arccosh}\left(1+\frac{2|\sigma_2-\sigma_1|^2}{(1-|\sigma_1|^2)(1-|\sigma_2|^2)}\right), \ \
\rm{for \ \  every} \ \ \sigma_1, \sigma_2\in \mathbb{B}^N.
\end{equation*}

For every $r>0$, we denote by
$
B(r)=\{\sigma\in \mathbb{B}^N\mid \, |\sigma|<r\}$
(resp. $B_H(r)=\{\sigma\in \mathbb{B}^N\mid\, d_H(\sigma)<r\})$
the Euclidean (resp. geodesic) ball of radius $r,$ at the origin $\sigma_0\in\mathbb{B}^N.$
It follows by \eqref{distance0}, that for every $r\in (0,1)$,
$
B(r)=B_H\left(\log\left(\frac{1+r}{1-r}\right)\right).$
See \cite{punzo} for additional comments and related facts.

\indent
For any $\sigma\in \R^N$, let $T_\sigma (\mathbb{B}^N)$ be the tangent space at $\sigma\in \mathbb{B}^N,$
equipped by the inner product $\left\langle \cdot, \cdot \right\rangle_\sigma$ and let $T(\mathbb{B}^N)=\bigcup_{\sigma\in \mathbb{B}^N}T_\sigma (\mathbb{B}^N)$ be the tangent bundle.
Whenever possible, given $X,Y\in T_\sigma(\mathbb{B}^N)$, we  write $|X|$ and $\langle X,Y\rangle$ instead of
$|X|_\sigma$ and  $g_\sigma(X,Y) = \langle X,Y\rangle_\sigma$, respectively.

Recall that $C^\infty_0(\mathbb{B}^N)$ denotes the space of real-valued smooth functions compactly supported on $\mathbb{B}^N$. Let
\begin{equation}\label{normag}
	\left\|u\right\|=\sqrt{\int_{\mathbb{B}^N}  |\nabla_H u(\sigma)|^2 d\mu}, \ \ \mbox{for every} \ \  u\in C^\infty_0(\mathbb{B}^N),
\end{equation}
 where $d\mu$ denotes the Riemannian measure on $\mathbb{B}^N$ from \eqref{dmu} and we get the following
$$
\nabla_H=\left(\frac{(1-|\sigma|^2)}{2}\right)^2\nabla
\quad\mbox{and}\quad |\nabla_H u(\sigma)|=\left(\frac{(1-|\sigma|^2)}{2}\right)^2\sqrt{\left\langle \nabla u(\sigma),\nabla u(\sigma)\right\rangle}.
$$
Then $H^{1,2}(\mathbb{B}^N)$ is the completion of $C_0^\infty(\mathbb{B}^N)$ with respect to the norm \eqref{normag} and it is a Hilbert space with the  inner product
\begin{equation}\label{prodscalareg}
	\left\langle u,v\right\rangle =
	\int_{\mathbb{B}^N} \left\langle \nabla_H u(\sigma),\nabla_H v(\sigma) \right\rangle d\mu, \ \
	\mbox{for every} \ \
	u,v\in H^{1,2}(\mathbb{B}^N).
\end{equation}

\indent We need to find critical points of the functional $J_\lambda$ from \eqref{Jlambda} so we shall invoke the principle of symmetric criticality, together with the following critical point theorem of Ricceri \cite{Ricceri} which we state in a form more suitable for our purpose.

\begin{theorem}\label{BMB}
Let $X$ be a reflexive real Banach space and  $\Phi,\Psi:X\to\R$
 G\^{a}teaux differentiable functionals such that  $\Phi$ is
strongly continuous, sequentially weakly lower semicontinuous and coercive,
whereas $\Psi$ is sequentially weakly upper semicontinuous.
 Given $r>\inf_X
\Phi$, let
$$
\varphi(r)=\inf_{u\in\Phi^{-1}((-\infty,r))}\frac{\displaystyle \sup\{\Psi(v)\mid {v\in\Phi^{-1}((-\infty,r))}  \}-\Psi(u)}{r-\Phi(u)}.
$$
Then for every $r>\inf_X\Phi$ and 
$\lambda\in\left(0,\frac{1}{\varphi(r)}\right)$, the restriction
of the functional $J_\lambda=\Phi-\lambda\Psi$ on
$\Phi^{-1}((-\infty,r))$ admits a global minimum which is a
critical point $($local minimum$)$ of $J_\lambda$ in $X$.
\end{theorem}

\begin{remark}
Theorem 2.1 is a direct consequence of \cite[Theorem 2.1]{BoMo} (see also 
 \cite{Bo,BoCa} for related topics).
 \end{remark}
\begin{remark}
Problem \eqref{BWP6} is set on the entire noncompact space $\mathbb{B}^N$. Therefore  we shall take
a group theoretical approach   in Section \ref{Section3}, in order to identify those symmetric subspaces of $H^{1,2}(\mathbb{B}^N)$ on which compactness of the embedding into $L^\nu(\mathbb{B}^N)$ can be regained.
\end{remark}
\section{$SO(N)$-invariant functions}\label{Section3}

Consider   the special orthogonal group $SO(N)$, $N\geq 3.$ Let $\cdot: SO(N)\times \mathbb{B}^N\rightarrow \mathbb{B}^N$ be the natural  action of $SO(N)$ on $\mathbb{B}^N$.
 The action $\circledast_{SO(N)}: SO(N)\times H^{1,2}(\mathbb{B}^N)\rightarrow H^{1,2}(\mathbb{B}^N)$ of a subgroup $SO(N)\in \mathscr F$ on $H^{1,2}(\mathbb{B}^N)$ is given by
 \begin{equation}\label{AsymmetricAction}
 g\circledast_{SO(N)}u(\sigma)=u(g^{-1}\cdot\sigma), \quad \mbox{ for a.e. } \;\sigma\in\mathbb{B}^N,
 \end{equation}
 for every $g\in SO(N)$ and $u\in H^{1,2}(\mathbb{B}^N)$.
Denote by
$$
H_{SO(N)}^{1,2}(\mathbb{B}^N)=\left\lbrace u\in H^{1,2}(\mathbb{B}^N)\mid g\circledast_{SO(N)}u=u, \mbox{ for every }g\in SO(N)\right\rbrace
$$
the subspace of $SO(N)$-invariant functions of $H^{1,2}(\mathbb{B}^N)$.
By using a recent embedding theorem of Skrzypczak and Tintarev \cite[Theorem 1.3 and Proposition 3.1]{SkTin}, the following compactness argument can be proved (see also  \cite{FaKr,Lions}).

\begin{theorem}\label{Compactness}{\rm{(See \cite{SkTin})}}
For every $\nu\in(2,2^*),$
the embedding $H_{SO(N)}^{1,2}(\mathbb{B}^N)\to L^\nu(\mathbb{B}^N)$ is compact.
\end{theorem}

Next, we recall the Palais principle of symmetric criticality. The group $(SO(N), * )$ acts continuously on the Hilbert space $H^{1,2}(\mathbb{B}^N)$ by  $(\tau,u)\mapsto \tau \circledast_{SO(N)} u$ from $SO(N)\times H^{1,2}(\mathbb{B}^N)$ to $H^{1,2}(\mathbb{B}^N)$, if this map itself is continuous on $SO(N)\times H^{1,2}(\mathbb{B}^N)$ and it has the following properties
 \begin{itemize} 	
 	\item[$(i_1)$] for every $\tau\in SO(N)$,
 	$u\mapsto \tau \circledast_{SO(N)} u$ is linear;
 	\item[$(i_2)$]  for every $\tau_1,\tau_2\in SO(N)$ and $u\in H^{1,2}(\mathbb{B}^N)$,
 	$$(\tau_1*\tau_2)\circledast_{SO(N)}u=\tau_1\circledast_{SO(N)}(\tau_2 \circledast_{SO(N)}u); \ \ \mbox{and}$$
 	\item[$(i_3)$]  for every $u\in H^{1,2}(\mathbb{B}^N)$, 
 	$$id_{SO(N)}\circledast_{SO(N)}u=u,$$
 	 where $id_{SO(N)}\in SO(N)$ denotes the identity element of $SO(N)$.
 \end{itemize}
 Define
 $$
\rm{Fix}_{SO(N)}(H^{1,2}(\mathbb{B}^N))=\left\lbrace u\in H^{1,2}(\mathbb{B}^N)\mid \tau\circledast_{SO(N)}u=u, \mbox{ for every }\tau\in SO(N)\right\rbrace
$$
and recall that the functional $\mathcal{J}:H^{1,2}(\mathbb{B}^N)\to\mathbb{R}$ is called $SO(N)$-invariant if
$$
\mathcal{J}(\tau\circledast_{SO(N)}u)=\mathcal{J}(u),
\ \ \mbox{for every} \ \ u\in H^{1,2}(\mathbb{B}^N).$$ The following result holds.

 \begin{theorem}\label{PalaisTh}{\rm{(See \cite{palais})}}
 	Let $H^{1,2}(\mathbb{B}^N)$ be the Sobolev space associated to the Poincaré model $\mathbb{B}^N$, $SO(N)$  the special orthogonal group acting continuously on $H^{1,2}(\mathbb{B}^N)$ by the map 
 	$$\circledast_{SO(N)}:SO(N)\times H^{1,2}(\mathbb{B}^N)\to H^{1,2}(\mathbb{B}^N),$$
 	 and $\mathcal{J}:H^{1,2}(\mathbb{B}^N)\to\mathbb{R}$  a $SO(N)$-invariant $C^1$-function. 
 	 
 	 If $u\in \rm{Fix}_{SO(N)}(H^{1,2}(\mathbb{B}^N))$ is a critical point of $\mathcal{J}_{|\rm{Fix}_{SO(N)}(H^{1,2}(\mathbb{B}^N))}$, then $u\in H^{1,2}(\mathbb{B}^N)$ is also a critical point of $\mathcal{J}$.
 	\end{theorem}
For details and comments we refer to
\cite[Section 5]{cp} and
\cite{demo}. See also \cite{Molica, MoPu1, puccirad} for related topics and results.

\section{Proof of the main theorem}\label{Section4}
Consider the functional
$\mathcal J_{\lambda}(u)= \Phi(u)-\lambda \Psi|_{SO(N)}(u),$ $u \in H_{SO(N)}^{1,2}(\mathbb{B}^N),$
where
$$\Phi(u)=\frac 1 2 \int_{\mathbb{B}^N}  |\nabla_H u(\sigma)|^2 d\mu \ \ \rm{and} \ \ \Psi(u)= \int_{\mathbb{B}^N} \alpha(\sigma)F(u(\sigma))d\mu.$$
We shall apply Theorem \ref{BMB} to
the energy functional~$\mathcal J_{\lambda}$ and use some ideas from  \cite{Molica,MoSec}.
On the basis of the preliminaries collected in Sections \ref{Section2} and \ref{Section3}, the existence of one
nontrivial $SO(N)$-symmetric solution of problem \eqref{BWP6} follows by the Palais criticality principle (Theorem \ref{PalaisTh}).

The space $H_{SO(N)}^{1,2}(\mathbb{B}^N)$ admits a Hilbert structure. By \cite{HoffSpru}, the functionals~$\Phi$ and
$\Psi|_{H_{SO(N)}^{1,2}(\mathbb{B}^N)}$ satisfy all the regularity assumptions of Theorem \ref{BMB}. More precisely, the functional $\Phi$ is (strongly) continuous, coercive in the symmetric space $H_{SO(N)}^{1,2}(\mathbb{B}^N)$ and
$\inf\{\Phi(u)\mid  u\in H_{SO(N)}^{1,2}(\mathbb{B}^N) \}=0.$

Since for every $\nu\in[2,2^*]$,  the Sobolev embedding $H^{1,2}(\mathbb{B}^N)\to L^\nu(\mathbb{B}^N)$ is continuous (but noncompact - see \cite{HoffSpru}), we shall make use also of the positive constant
$$
c_\nu=\sup
\Big\{
{\left( \int_{\mathbb{B}^N} |u(\sigma)|^\nu  d\mu \right)^{\frac{1}{\nu}}}
{\left( \int_{\mathbb{B}^N}  |\nabla_H u(\sigma)|^2 d\mu \right)^{-\frac{1}{2}}}
\mid
u\in H^{1,2}(\mathbb{B}^N)\setminus\{0\}
\Big\}.
$$

Set
$
h(\omega)= {\omega} (q \sqrt{2} \| \alpha \|_{p}+2^{q/2}c_{q}^{q-1} \| \alpha \|_{\infty}{\omega}^{q-1})^{-1},
$
for every $\omega>0$,
 and define
\begin{equation}\label{lambda1}
\lambda^{\star}=\frac{q\max\{h(\omega)\mid  \omega>0 \}}{\alpha_fc_q}, \ \ \mbox{where} \ \  p=\frac{q}{q-1}.
\end{equation}
Take $0<\lambda<\lambda^{\star}$.
By \eqref{lambda1}, we have
\begin{equation}\label{n3}
\lambda<\lambda^{\star}{(\bar{\omega})}=\frac{qh(\bar \omega)}{\alpha_fc_q}, \ \ \mbox{for some} \ \  \bar{\omega}>0.
\end{equation}
\noindent Let $\Theta: (0,\infty)\rightarrow [0,\infty)$ be the real function defined by
$$
\Theta(r)=\frac{1}{r}  \sup\{\Psi|_{H_{SO(N)}^{1,2}(\mathbb{B}^N)}(u)\mid u\in\Phi^{-1}((-\infty,r)) \}, \ \ \mbox{for every} \ \  r>0.
$$
 Then condition \eqref{Growth} gives
$$
\Psi(u)\leq \alpha_f\int_{{\mathbb{B}^N}}\alpha(\sigma)|u(\sigma)|d\mu+\frac{\alpha_f}{q}\int_{{\mathbb{B}^N}}\alpha(\sigma)|u(\sigma)|^qd\mu, \ \ \mbox{for every} \ \  u\in {H_{SO(N)}^{1,2}(\mathbb{B}^N)},
$$
 so if $\Phi(u)<r$, then
\begin{equation}\label{gio}
  \int_{\mathbb{B}^N}  |\nabla_H u(\sigma)|^2 d\mu< {2r}, \ \ \mbox{ for every} \ \ u\in H_{SO(N)}^{1,2}(\mathbb{B}^N).
\end{equation}
An application of (\ref{gio}) and Theorem \ref{Compactness} yields
$$
\int_{\mathbb{B}^N} \alpha(\sigma)F(u(\sigma))d\mu< \alpha_fc_q\left(\|\alpha\|_p\sqrt{{2r}}+\frac{c_q^{q-1}}{q}\|\alpha\|_{\infty}(2r)^{q/2}\right), \ \ \mbox{for every} \ \  u\in H_{SO(N)}^{1,2}(\mathbb{B}^N), \Phi(u)<r.
$$
 Consequently,
 $$
\sup\{\Psi|_{H_{SO(N)}^{1,2}(\mathbb{B}^N)}(u)\mid u\in\Phi^{-1}((-\infty,r)) \}
\leq
\alpha_fc_q\left(\|\alpha\|_p\sqrt{{2r}}+\frac{c_q^{q-1}}{q}\|\alpha\|_{\infty}(2r)^{q/2}\right).
$$
Hence
\begin{equation}\label{n}
\Theta(r)\leq \alpha_fc_q\left(\|\alpha\|_p\sqrt{\frac{2}{r}}+\frac{2^{q/2}c_q^{q-1}}{q}\|\alpha\|_{\infty}r^{q/2-1}\right),
\ \ \mbox{ for every} \ \  r>0.
\end{equation}
Taking $r=\bar\omega^2$, we get
\begin{equation}\label{n2}
\Theta(\bar{\omega}^2)\leq \alpha_fc_q\left(\sqrt{{2}}\frac{\|\alpha\|_p}{\bar\omega}+\frac{2^{q/2}c_q^{q-1}}{q}\|\alpha\|_{\infty}{\bar\omega}^{q-2}\right).
\end{equation}
On the other hand,
$$
\varphi(\bar{\omega}^2)
=
\inf_{u\in\Phi^{-1}((-\infty,\bar{\omega}^2))}
\frac{\displaystyle\sup\{\Psi|_{H_{SO(N)}^{1,2}(\mathbb{B}^N)}(u)\mid u\in\Phi^{-1}((-\infty,r)) \}
-\Psi|_{H_{SO(N)}^{1,2}(\mathbb{B}^N)}(u)}
{r-\Phi(u)}
\leq
 \Theta(\bar{\omega}^2),
$$
since $0\in\Phi^{-1}((-\infty,\bar{\omega}^2))$ and $\Phi(0)=0$.
 By virtue of (\ref{n3}) and \eqref{n2}, we have
\begin{equation}\label{n2fg}
\varphi(\bar{\omega}^2)\leq \Theta(\bar{\omega}^2)<\frac{1}{\lambda},
\end{equation}
hence
$
\lambda\in (0,{1}/{\varphi(\bar{\omega}^2)}).
$
Consequently, by Theorem \ref{BMB}, there exists  $u_\lambda^{SO(N)}\in\Phi^{-1}((-\infty,\bar{\omega}^2))$ such that
$$
\Phi'(u_\lambda^{SO(N)})=\lambda (\Psi|_{H_{SO(N)}^{1,2}(\mathbb{B}^N)})'(u_\lambda^{SO(N)}).
$$
Moreover, $u_\lambda^{SO(N)}$ is a global minimum of $\mathcal{J}_{\lambda}$ on the sublevel $\Phi^{-1}((-\infty,\bar{\omega}^2))$.

Next, we  show that
solution~$u_\lambda^{SO(N)}$  is not the trivial (identically zero) function. If
$f(0)\not =0$, then it easily follows that $u_\lambda^{SO(N)}\not \equiv
0$ in $H_{SO(N)}^{1,2}(\mathbb{B}^N)$, since the trivial function does not solve
problem~\eqref{BWP6}. 

So let us consider the case when $f(0)=0$ and fix
$\lambda\in (0, \lambda^{\star}{(\bar{\omega})})$ for some $\bar \omega>0$. Let $u_\lambda^{SO(N)}$ be such that
\begin{equation}\label{minimum}
\mathcal J_{\lambda}(u_\lambda^{SO(N)})\leq \mathcal
J_{\lambda}(u),\quad \mbox{for every}\,\,\, u\in H_{SO(N)}^{1,2}(\mathbb{B}^N)\,\,
\mbox{such that}\,\,\,\Phi(u)<\bar{\omega}^2
\end{equation}
and
\begin{equation}\label{minimum2}
\Phi(u_\lambda^{SO(N)})<\bar\omega^2\,,
\end{equation}
and that $u_\lambda^{SO(N)}$ is  a critical point of $\mathcal
J_{\lambda}$ in $H_{SO(N)}^{1,2}(\mathbb{B}^N)$.

 Applying Theorem \ref{MainTheorem}, the energy $J_\lambda$ defined in \eqref{Jlambda} needs to be invariant with respect to the special orthogonal group $SO(N)$. To show this,  fix $u\in H^{1,2}(\mathbb{B}^N)$ and $g\in SO(N)$. Since $g\in SO(N)$ is an isometry, it follows by \eqref{AsymmetricAction} that
\begin{equation}\label{chain}
\nabla_H (g\circledast_{SO(N)}u)(\sigma)=Dg_{g^{-1}\cdot\sigma}\nabla_H u(g^{-1}\cdot\sigma),
\ \ \mbox{for a.e.} \ \ \sigma \in \mathbb{B}^N.
\end{equation}
 If $z=g^{-1}\cdot \sigma$, then
\begin{align}\label{inv1}
	\left\|g\circledast_{SO(N)}u\right\|^2 & =\int_{\mathbb{B}^N}  |\nabla_H (g\circledast_{SO(N)}u)(\sigma)|_\sigma^2   d\mu(\sigma)\nonumber\\
	&
	 = \int_{\mathbb{B}^N}  |\nabla_H u(g^{-1}\cdot\sigma)|_{g^{-1}\cdot\sigma}^2   d\mu(\sigma)
 = \left\|u \right\|^2.
	\end{align}
where we have  used \eqref{chain}.
On the other hand, since $\alpha \in L^1(\mathbb{B}^N)\cap L^{\infty}(\mathbb{B}^N)$ is radially symmetric respect to the origin, it follows that
\begin{align}\label{inv2}
\int_{\mathbb{B}^N} \alpha(\sigma)\left(\int_0^{(g\circledast_{SO(N)}u)(\sigma)}h(t)dt\right)d\mu(\sigma)	& = \int_{\mathbb{B}^N} \alpha(\sigma)\left(\int_0^{u(g^{-1}\cdot\sigma)}h(t)dt\right)d\mu(\sigma)\\
	& =\int_{\mathbb{B}^N} \alpha(z)\left(\int_0^{u(z)}h(t)dt\right)d\mu(z).\nonumber
	\end{align}
By \eqref{inv1} and \eqref{inv2},  we have
$
J_{\lambda}(g\circledast_{SO(N)}u)=J_{\lambda}(u),
$
which proves the $SO(N)$ invariance of the functional $J_{\lambda}$.

By Theorem \ref{PalaisTh}, it is clear that $u_\lambda^{SO(N)}$ weakly solves
problem~\eqref{BWP6}\,.
Proving that $u_\lambda^{SO(N)}\not \equiv 0$ in $H_{SO(N)}^{1,2}(\mathbb{B}^N)$\,,
we show the existence of a sequence $\big\{w_j\big\}_{j\in
\NN}$ in $H_{SO(N)}^{1,2}(\mathbb{B}^N)$ such that
\begin{equation}\label{wjBlu}
\limsup_{j\to
\infty}\frac{\Psi|_{H_{SO(N)}^{1,2}(\mathbb{B}^N)}(w_j)}{\Phi(w_j)}=\infty\,.
\end{equation}
By (\ref{azero}), there exists
$\{t_j\}_{j\in\NN}\subset (0,+\infty)$ such that $t_j\to 0^+$ 
when
 $j\to
+\infty$, and
\begin{equation}\label{limsup}
 \lim_{j\rightarrow +\infty}\frac{\displaystyle F(t_j)}{t_j^2}=+\infty.
\end{equation}
Therefore
for every $M>0$ and
all
sufficiently large  $j$,
\begin{equation}\label{M}
\displaystyle\,F(t_j)>Mt_j^2\,.
\end{equation}
Now, $\alpha\in L^{\infty}(\mathbb{B}^N)\setminus\{0\}$ is nonnegative in $\mathbb{B}^N$. Hence there are real numbers $\rho>r>0$ and $\alpha_0>0$ such that
\begin{equation}\label{ess}
\textrm{essinf}_{\sigma\in A_{r}^{\rho}}\alpha(\sigma)\geq \alpha_0>0.
\end{equation}
For every $0<a<b$, set
\[
A_a^{b}=\left\{\sigma\in \mathbb{B}^N\mid b-a< \log\left(\frac{1+|\sigma|}{1-|\sigma|}\right)< a+b \right\}.
\]
Define $w^{1/2}_{\rho,r}(\sigma)\in H^{1,2}(\mathbb{B}^N)$  by
\begin{equation} \label{TestFunction}
w^{1/2}_{\rho,r}(\sigma)=\left\{
\begin{array}{ll}
1 & \mbox{ if $\sigma \in A_{1/2 r}^{\rho}$}\\[5pt]
0 & \mbox{ if $\sigma \in \mathbb{B}^N \setminus A_r^{\rho}$} \\[5pt]
 \displaystyle\frac{2}{r}\left({r- \left|\displaystyle\log\left(\frac{1+|\sigma|}{1-|\sigma|}\right)-\rho\right|}\right)
& \mbox{ if $\sigma \in A_{r}^{\rho}\setminus A_{1/2 r}^{\rho}$},
\end{array}
\right.
\end{equation}
\noindent for every $\sigma\in \mathbb{B}^N$.

Since the group $SO(N)$ is a compact connected subgroup of the isometry group $\textrm{Isom}_g(\mathbb{B}^N)$ such that $\Fix_{SO(N)}(\mathbb{B}^N)=\{\sigma_0\}$, it follows  that $w^{1/2}_{\rho,r}\in H^{1,2}(\mathbb{B}^N)$, given in \eqref{TestFunction}, belongs to $H^{1,2}_{SO(N)}(\mathbb{B}^N)$. Therefore
  $\operatorname{supp}(w^{1/2}_{\rho,r})\subseteq A_r^{\rho}(\sigma_0)$,
  $\|w^{1/2}_{\rho,r}\|_\infty\leq 1$, and
 $w^{1/2}_{\rho,r}(\sigma)=1$, for every $\sigma\in A_{1/2 r}^{\rho}(\sigma_0)$. 

\begin{remark}
The test functions used here were introduced in \cite{MolicaNA}, following \cite{FaKr}. We note that test functions introduced in \cite{Bon,BoLi} are different. We also emphasize that the different geometrical structure used along the proof is crucial in order to recover the $SO(N)$ invariance of the test functions.
\end{remark}
 Define $w_j=t_jw^{1/2}_{\rho,r}$ for any $j\in \NN$\,. Taking into account that $w^{1/2}_{\rho,r}\in H_{SO(N)}^{1,2}(\mathbb{B}^N),$ it is easily seen that $w_j\in H_{SO(N)}^{1,2}(\mathbb{B}^N)$, for every $j\in \NN$. Furthermore, exploiting the properties of $w^{1/2}_{\rho,r}$, by (\ref{M}), it follows that:
\begin{eqnarray}\label{conto1}
& {\frac{\displaystyle \Psi|_{H_{SO(N)}^{1,2}(\mathbb{B}^N)}(w_j)}{\displaystyle\Phi(w_j )}}= \frac{\displaystyle \int_{A_{1/2 r}^{\rho}}\alpha(\sigma)F(w_j(x))\, d\mu+\int_{A_{r}^{\rho}\setminus A_{1/2 r}^{\rho}} \alpha(\sigma)F(w_j(\sigma))\,d\mu}{\displaystyle \Phi(w_j)} \nonumber \\
& \qquad \qquad= \frac{\displaystyle \int_{A_{1/2 r}^{\rho}}\alpha(\sigma)F(t_j)\, d\mu+\int_{A_{r}^{\rho}\setminus A_{1/2 r}^{\rho}}
\alpha(\sigma)F(t_jw^{1/2}_{\rho,r}(\sigma))\,d\mu}{\displaystyle\Phi(w_j)}  \\
& \qquad \qquad\quad\geq \displaystyle 2\alpha_0\frac{M\mu(A_{1/2 r}^{\rho})t^2_j+{\displaystyle \int_{A_{r}^{\rho}\setminus A_{1/2 r}^{\rho}}
F(t_jw^{1/2}_{\rho,r}(\sigma))\,d\mu}}{t_j^2\|w^{1/2}_{\rho,r}\|^2}, 
\ \ \mbox{for  sufficiently large} \ \  j. \nonumber
\end{eqnarray}
Assertion~(\ref{wjBlu}) now follows by (\ref{conto1}).

Now
$$\int_{\mathbb{B}^N}  |\nabla_H w_j(\sigma)|^2 d\mu=t_j^{2}\,\int_{\mathbb{B}^N}  |\nabla_H w^{1/2}_{\rho,r}(\sigma)|^2 d\mu\to 0,
\ \ \mbox{as} \ \  j\to +\infty,
$$
 so
$\|w_j\|< \sqrt{2}\bar \omega,$
for sufficiently large $j$. Hence
\begin{equation}\label{wj1}
w_j\in \Phi^{-1}\big((-\infty, \bar \omega^2)\big)\,,
\end{equation}
provided that $j$ is large enough. Moreover, by \eqref{wjBlu},
\begin{equation}\label{wj2}
\mathcal
J_{\lambda}(w_j)=\dfrac{1}{2}\int_{\mathbb{B}^N}  |\nabla_H w_j(\sigma)|^2 d\mu-\lambda\int_{\mathbb{B}^N} \alpha(\sigma)\left(\int_0^{w_j(\sigma)}f(t)dt\right)d\mu<0,
\end{equation}
for sufficiently large $j$ and $\lambda>0$.

Since the restriction of $\mathcal
J_{\lambda}$ to $\Phi^{-1}\big((-\infty,\bar \omega^2)\big)$ has $u_\lambda^{SO(N)}$ as a global minimum, it follows by (\ref{wj1}) and (\ref{wj2}) that
\begin{equation}\label{nontrivial}
\mathcal J_{\lambda}(u_\lambda^{SO(N)})\leq \dfrac{1}{2}\int_{\mathbb{B}^N}  |\nabla_H w_j(\sigma)|^2 d\mu-\lambda\int_{\mathbb{B}^N} \alpha(\sigma)\left(\int_0^{w_j(\sigma)}f(t)dt\right)d\mu<\mathcal J_{\lambda}(0)\,,
\end{equation}
so  $u_\lambda^{SO(N)}\not\equiv 0$ in $H_{SO(N)}^{1,2}(\mathbb{B}^N)$ as asserted.

 Therefore $u_\lambda^{SO(N)}$ is a
nontrivial weak solution of problem~\eqref{BWP6}. The arbitrariness
of $\lambda$ implies that $u_\lambda^{SO(N)}\not \equiv 0$, for every
$\lambda\in (0, \lambda^{\star})$.

Finally, we show that $\displaystyle \lim_{\lambda\rightarrow
0^+}{\|u_\lambda^{SO(N)}\|}=0.$
To this end, consider
$\lambda\in (0, \lambda^{\star}{(\bar{\omega})})$ for some $\bar \omega>0$. Taking into account that $\Phi(u_\lambda^{SO(N)})<\bar \omega^2$, it follows that
$\Phi(u_\lambda^{SO(N)})=\frac 1 2\|u_\lambda^{SO(N)}\|^2<\bar
\omega^2,$ i.e.,
$\|u_\lambda^{SO(N)}\|<\sqrt{2}\bar \omega.$

The growth condition (\ref{Growth}) yields
\begin{eqnarray*}\label{bounded}
& {\displaystyle \left|\Psi'(u_\lambda^{SO(N)})\right|  \leq\alpha_f\Bigg(
\int_{\mathbb{B}^N}\alpha(\sigma)|u_\lambda^{SO(N)}(\sigma)|d\mu+\int_{\mathbb{B}^N}\alpha(\sigma)|u_\lambda^{SO(N)}(\sigma)|^qd\mu}\Bigg) \nonumber\\
& \qquad \qquad\qquad {\displaystyle \leq \alpha_f\Bigg(\|\alpha\|_p\|u_\lambda^{SO(N)}\|_q+\|\alpha\|_\infty\|u_\lambda^{SO(N)}\|_q^q}\Bigg) \\
& \qquad \qquad \qquad\quad \qquad\qquad\,{\displaystyle < c_q\alpha_f\Bigg(\sqrt{2}\|\alpha\|_p\bar \omega+2^{q/2}c_q^{q-1}\|\alpha\|_\infty\bar\omega^q}\Bigg)=M_{\bar \omega}\,.\nonumber
\end{eqnarray*}
Since $u_\lambda^{SO(N)}$ is a critical point of $\mathcal
J_{\lambda}$\,, it follows that $\langle \mathcal
J_{\lambda}'(u_\lambda^{SO(N)}), \varphi \rangle=0$, for every $\varphi
\in H_{SO(N)}^{1,2}(\mathbb{B}^N)$ and  $\lambda\in (0,\lambda^{\star}{(\bar{\omega})})$. Hence, $\langle
\mathcal J_{\lambda}'(u_\lambda^{SO(N)}), u_\lambda^{SO(N)}\rangle=0$
and thus
\begin{equation*}\label{rangle}
\langle \Phi'(u_\lambda^{SO(N)}), u_\lambda^{SO(N)}\rangle=\lambda \Psi'(u_\lambda^{SO(N)}), \ \ \rm{for \ \ every} \ \
\lambda\in (0,\lambda^{\star}{(\bar{\omega})}).
\end{equation*}
The  relations above now ensure that
$$\displaystyle 0\leq \|u_\lambda^{SO(N)}\|^2=\langle\Phi'(u_\lambda^{SO(N)}),u_\lambda^{SO(N)}\rangle=\lambda\Psi'(u_\lambda^{SO(N)})< \lambda M_{\bar \omega}, \ \ \rm{for \ \ every} \ \
\lambda\in (0,\lambda^{\star}{(\bar{\omega})}).$$ 
Hence $\displaystyle
\lim_{\lambda\rightarrow 0^+}{\|u_\lambda^{SO(N)}\|}=0$, as was
asserted.
The proof of Theorem~\ref{MainTheorem} is now complete.\hfill $\Box$

\begin{remark}
Profile decomposition methods can be useful in order to study similar problems when a lack of compactness occurs (see, among others, the recent papers \cite{De2,De1}).
A further and more general investigation of this topics will be included in the forthcoming book \cite{MolicaPucciBook}. 
\end{remark}

\begin{remark}
The referee has observed that the importance of the solution as a local minimum is in that we can obtained in addition a second solution, and suggested as a further study, to attempt to apply the result contained in \cite{BoDa} to obtain two nonzero solutions for this type of problems. 
\end{remark}

\section{An example}\label{Section5}
 We conclude the paper by exhibiting the following model equation which illustrates how our main result can be applied.

\begin{ex} 
For any $1<r<2$, consider the following problem
on ${\mathbb{B}^4}$
\begin{equation}\label{eq}
-\Delta_{H}u=\lambda \left(\frac{1-|\sigma|^2}{2}\right)^4  |u|^{r-2}u, \ \ u\in H^{1,2}(\mathbb{B}^4).
\end{equation}
By Theorem \ref{MainTheorem}, there exists $\lambda^{\star}>0$ such that for every $\lambda\in (0,\lambda^{\star}_{SO(N)})$, problem~\ref{eq} 
admits at least one nontrivial $SO(N)$-symmetric weak solution
$
u_\lambda^{SO(N)}\in H^{1,2}(\mathbb{B}^4) $
such  that $
\lim_{\lambda\rightarrow 0^+}\|u_\lambda^{SO(N)}\|=0.
$
\end{ex}

\medskip

\section*{Acknowledgments}
This work was stimulated by discussions with G. Molica Bisci whom we thank for comments and suggestions. The research of the first and the fourth author  reported in this paper was supported by the Slovenian Research Agency  (P1-0292, N1-0114, N1-0083, N1-0064, J1-8131).
The research of the second and the third author  reported in this paper, carried out at the Budapest University of Technology and Economics, was supported by the "National Challenges Program" of the National Research Development and Innovation Office (BME NC TKP2020) and the Hungarian Ministry of Human Capacities OTKA (SNN125119).
We thank the referee for comments and suggestions.

\end{document}